\documentclass[12pt,oneside]{amsart}
\usepackage{amsmath, amssymb, verbatim}

\theoremstyle{plain}
\newtheorem{lemma}{Lemma}

\newtheorem*{maintheorem}{Theorem}

\begin{document}

\title[DIFFERENCE SETS AND 	SHIFTED PRIMES]{DIFFERENCE SETS AND 	SHIFTED PRIMES}
\author{Jason Lucier}
\date{May 24, 2007}

\maketitle

\section{Introduction}\label{intro}
For a set of integers $A$ we denote by $A-A$ the set of all differences $a-a'$ with $a$ and $a'$ in $A$,
and if $A$ is a finite set we denote its cardinality by $|A|$.  
S\'{a}rk\"{o}zy \cite{Sarkozy1} proved, by the Hardy-Littlewood method, 
that if $A$ is a subset of $\{1, \ldots, n\}$ such that $A-A$ does not contain a 
perfect square, then
\[
 	|A| \ll n(\log_2 n)^{2/3}(\log n)^{-1/3}.
\]
This estimate was improved by Pintz, Steiger and Szemer\'{e}di \cite{PSS} to 
\[
 	|A| \ll n(\log n)^{-(1/12)\log\log\log\log n}.
\]
This improvement was obtained using the Hardy-Littlewood method together 
with a combinatorial result concerning sums of rationals.
Balog, Pelik\'{a}n, Pintz and Szemer\'{e}di \cite{BPPS}, elucidating the method in \cite{PSS}, 
proved for any fixed integer $k \geq 2$, that if $A$ is a subset of $\{1, \ldots , n\}$ such 
that $A-A$ does not contain a perfect $k$-th power, then
\[
  |A| \ll_k n(\log n)^{-(1/4)\log\log\log\log n}.
\]

In the works cited above the following basic property is used;
if $s$ is a perfect $k$-th power then so is $q^{k}s$ for every positive integer $q$. This multiplicative 
property is used in the following fashion:
Suppose that $B$ is a set of integers and 
$A = \{c+q^k b : b \in B \}$ for some integers $c$ and $q \geq 1$, 
if $A-A$ does not contain a perfect $k$-th power, then the same is true for $B-B$. 
This deduction is the basis of an iteration argument that plays a fundamental r\^{o}le in \cite{BPPS}, \cite{PSS}, and 
\cite{Sarkozy1}.

S\'{a}rk\"{o}zy \cite{Sarkozy3} also considered the set $\mathcal{S}= \{ \, p-1 \, : \, \text{$p$ a prime} \, \}$ of shifted primes, 
and showed that if $A$ is a subset of $\{1, \ldots , n\}$ such that $A-A$ does not contain an integer from $\mathcal{S}$ then 
\[
 		|A| \ll n \frac{(\log \log \log  n)^3(\log \log\log \log  n)}{(\log \log  n)^{2}}. 
\]
The argument S\'{a}rk\"{o}zy used in \cite{Sarkozy1} cannot be applied directly to the set $\mathcal{S}$ of shifted primes since 
it does not have a multiplicative property analogous to the one possessed by the set of perfect $k$-th powers. 
S\'{a}rk\"{o}zy got around this difficulty by not only considering the set $\mathcal{S}$ 
of shifted primes, but also the sets defined for each positive integer $d$ by
\[
	\mathcal{S}_d = \left\{ \; \frac{p-1}{d} \; : \; \text{$p$ a prime, $p \equiv 1 \pmod{d}$} \; \right\}.
\]
In \cite{Sarkozy3} S\'{a}rk\"{o}zy uses an iteration argument based on the following observation.
Suppose $B$ is a set of integers and $A = \{c+q b : b \in B \}$ for some integers $c$ and $q \geq 1$, 
if $A-A$ does not intersect $S_d$ for some positive integer $d$, then $B-B$ does not intersect $\mathcal{S}_{dq}$. 

In this article we show that the combinatorial argument presented in \cite{BPPS} and \cite{PSS} can be
carried out to improve S\'{a}rk\"{o}zy's result on the set $\mathcal{S}$ of shifted primes.
We shall prove the following.

\begin{maintheorem}
Let $n$ be a positive integer and $A$ a subset of $\{1, \ldots , n\}$. If there does not exist a pair of  
integers $a,a' \in A$ such that $a-a'=p-1$ for some prime $p$, then 
\[
|A| \ll 
 	n\left(\frac{(\log\log\log n)^3(\log\log\log\log n)}{(\log\log n)} \right)^{\log\log\log\log\log n}.
\]
\end{maintheorem}

The set of perfect squares and the set $\mathcal{S}$ of shifted primes are examples of \textit{intersective} sets.
To define this class of sets we introduce some notation. Given a set of positive integers $H$ we define $D(H,n)$, 
for any positive integer $n$, to be the maximal size of a subset $A$ of $\{1, \ldots , n\}$ such that $A-A$ does 
not intersect $H$. A set of positive integers $H$ is called \textit{intersective} if $D(H,n) = o(n)$.

Kamae and Mend\`{e}s France \cite{Kamae_&_MendesFrance} supplied a general criterion for 
determining if a set of positive integers is intersective. From their criterion they deduced the following.
\begin{enumerate}
\item[(I)] For any fixed integer $a$ the set $\{\, p+a \, : \text{$p$ a prime}, p > -a \}$ 
					 is intersective if and only if $a = \pm 1$. 

\item[(II)]  Let $h$ be a nonconstant polynomial with integer coefficients and whose leading coefficient is positive. 
The set $\{\, h(m)\, : \, m \geq 1, h(m) \geq 1 \,\}$ is intersective if and only if for each positive integer $d$ the 
modular equation $h(x) \equiv 0 \pmod{d}$ has a solution.
\end{enumerate}

Let $h$ be a polynomial as in (II) with degree $k \geq 2$ and such that $h(x) \equiv 0 \pmod{d}$ has a solution for every positive integer $d$. The author \cite{Lucier} has shown that if $A$ is a subset of $\{1, \ldots ,n\}$ 
such that $A-A$ does not intersect $\{\, h(m)\, : \, m \geq 1, h(m) \geq 1 \,\}$, then 
$|A| \ll n(\log_2 n)^{\mu/(k-1)}(\log n)^{-(k-1)}$, where $\mu = 3$ if $k=2$ and $\mu =2$ if $k \geq 3$. 
It is possible to improve this result with the method presented in this paper.

\section{Preliminary lemmata}
In this paper we use the following notations. For a real number $x$ we write $e(x)$ for $e^{2\pi i x}$, and $[x]$ is 
used to denote the greatest integer less than or equal to $x$. The greatest common divisor of the integers $u$ and $v$ 
is given by $(u,v)$. Euler's totient function is given, as usual, by $\phi$. For any positive integer $i$ we write 
$\log_i$ to denote the $i$-th iterated logarithm, that is, $\log_1 n = \log n$ and $\log_i n = \log (\log _{i-1} n)$ 
for every integer $i \geq 2$. 

A fundamental r\^{o}le is played by the following relations; for integers $n$ and $r$, with $n$ positive,
\[
  \sum_{t=0}^{n-1} e(rt/n) = \begin{cases} n & \text{if $n|r$} \\
                              		   0 & \text{if $n \nmid r$}
                             \end{cases},
\quad \quad
  \int_0^1 e(r\alpha)d\alpha = \begin{cases} 1 & \text{if $r=0$} \\
                              		     0 & \text{if $r \ne 0$}
                             \end{cases}.
\]
Given a subset $A$ of $\{1, \ldots , n\}$ its generating function is given by 
\[
 	F(\alpha) = \sum_{a \in A} e(\alpha a), \quad \alpha \in \mathbb{R}.
\]
Using the relations above we find that 
\[
	\sum_{t = 1}^n |F(t/n)|^2 = n|A|, \quad \quad \int_0^1 |F(\alpha)|^2 d\alpha = |A|.
\]
Of course, these are particular cases of Parseval's identity.

S\'{a}rk\"{o}zy's method in \cite{Sarkozy1} and \cite{Sarkozy3} is based on Roth's work \cite{Roth} 
on three-term arithmetic progressions in dense sets. 
Following this method S\'{a}rk\"{o}zy uses a functional inequality to derive his results concerning the set of 
perfect squares and the set $\mathcal{S}$ of shifted primes. 
Our approach here uses, like Gowers \cite{Gowers} and Green \cite{Green}, 
a density increment argument. The next lemma tells us that if the generating function of a finite set 
$A$ satisfies a certain size constraint, then it must be concentrated along an arithmetic progression. 
We use this result in Lemma~\ref{crank} to obtain a density increment that we iterate in the final section 
of the paper to prove the theorem.
 
\begin{lemma}\label{bump} 
Let $n$ be a positive integer and $A$ a subset of $\{1, \ldots , n\}$ with size $\delta n$. For any real $\alpha$ let 
$F(\alpha)$ denote the generating function of $A$. 
Let $q$ be a positive integer and $U$ a positive real number such that $2\pi q U \leq n$. Let $E$ 
denote the subset of $[0,1]$ defined by
\[
E = \left\{ \alpha \in [0,1] : 
	\left|\alpha - \frac{a}{q} \right| \leq \frac{U}{n} \text{ for some $0 \leq a \leq q$} \right\}.
\]
If $\theta$ is a positive number such that 
\begin{equation}\label{tonka}
	\sum_{\substack{ t = 1 \\ t/n \in E}}^{n-1} \left|F(t/n)\right|^2 \geq \theta |A|^2,
\end{equation}
then there exists an arithmetic progression $P$ in $\{1, \ldots , n\}$ with difference $q$ such that 
\[
|P| \geq \frac{n}{32\pi q U} \quad \; \text{and} \; \quad |A \cap P| \geq  |P| \delta\big(1  + 8^{-1}\theta \big).
\]
\end{lemma}

\begin{proof}
This closely resembles Lemma 20 in \cite{Lucier} and can be proved in the same manner. 
\end{proof}

We now state a combinatorial result presented by Balog, Pelik\'{a}n, Pintz and Szemer\'{e}di in \cite{BPPS},
the proof of which uses only elementary techniques. It is this result, that we use in 
Lemma~\ref{lemma_nonuniform}, that allows us to improve S\'{a}rk\"{o}zy result on the set $\mathcal{S}$ of shifted primes.

\begin{lemma}\label{lemma_CR}
Let $K$ and $L$ be positive integers, and 
let $\tau$ be the maximal value of the divisor function up to $KL$. 
Let $\mathcal{K}$ be a nonempty subset of rationals such that if $a/k \in \mathcal{K}$ is in lowest terms then 
$1 \leq a \leq k \leq K$.  Suppose that for each $a/k \in \mathcal{K}$ there 
corresponds a subset of rationals $\mathcal{L}_{a/k}$ such that if $b/l \in \mathcal{L}_{a/k}$ is in lowest terms then 
$1 \leq b \leq l \leq L$. Suppose further that $B$ and $H$ are positive integers   
such that 
\[
 				|\mathcal{L}_{a/k}| \geq H \quad \text{for all $a/k \in \mathcal{K}$}
\]
and 
\[
 \left|\left\{ b : \frac{b}{l} \in \bigcup \mathcal{L}_{a/k} \right\}\right| \leq B \quad \text{for all $l \leq L$.}
\]
Then the size of the set 
\[
	\mathcal{Q} = \left\{ \; \frac{a}{k}+\frac{b}{l} \; : \; \frac{a}{k} \in \mathcal{K}, \, \frac{b}{l} \in \mathcal{L}_{a/k} \; \right\}
\]
satisfies
\[
 |\mathcal{Q}| \geq |\mathcal{K}| H\left(\frac{H}{LB\tau^8(1+\log K)}\right).
\]
\end{lemma}

\begin{proof}
This is Lemma CR in \cite{BPPS}.
\end{proof}

\section{Exponential sums over primes}
Let $d$ and $n$ denote positive integers. As in \cite{Sarkozy3}, our application of the Hardy-Littlewood method employs
exponential sums over numbers from the set $\mathcal{S}_d$ defined in the introduction. 
For any real number $\alpha$ we set 
\[
	S_{n,d}(\alpha) = \sum_{\substack{ s \in \mathcal{S}_d \\ s \leq n}} \log (ds+1) e(\alpha s).
\] 
In this section we present some estimates related to $S_{n,d}(\alpha)$. Throughout this section we assume $d$ and $n$ satisfy
\[
 d \leq \log n.
\]

\begin{lemma}\label{Siegel_Walfisz} 
For $n$ sufficiently large,
\[
		S_{d,n}(0) \gg \frac{dn}{\phi(d)}.
\]
\end{lemma}

\begin{proof}
By the definition of $\mathcal{S}_d$ we find that
\[
S_{d,n}(0) = \sum_{\substack{p \leq dn+1 \\ p \equiv 1 \mod d}} \log p.
\]
Since $d \leq \log n$ the Siegel-Walfisz theorem says that this sum is asymptotic to $(dn+1)/\phi(q)$, from which the result follows.
\end{proof}

The next two lemmas provide estimates of $S(\alpha)$ derived by A. S\'{a}rk\"{o}zy.

\begin{lemma}\label{majorestimates} 
Let $a$ and $b$ be integers such that $(a,b)=1$ and $1 \leq b \leq \log n$.
There exists a positive real number $c$ such that if  $\alpha$ is a real number that satisfies 
\begin{equation*}
	\left|\alpha - \frac{a}{b}\right| \leq \frac{\exp(c(\log n)^{1/2})}{n},
\end{equation*}
and $n$ is sufficiently large, then
\begin{equation*}
\left| S_{d,n}\left( \alpha \right) \right| <  \frac{dn}{\phi(d)\phi(b)},
\end{equation*}
furthermore, if $\alpha \ne a/b$ then 
\begin{equation*}
	\left| S_{d,n}\left( \alpha \right) \right| < \frac{d}{\phi(d)\phi(b)}\left|\alpha - \frac{a}{b}\right|^{-1}.  
\end{equation*}
\end{lemma}

\begin{proof}
This is a restatement of Lemma 5 from \cite{Sarkozy3}.
\end{proof}

Let $R$ denote a real number that satisfies 
\begin{equation}
 	3 \leq R \leq \log n.
\end{equation}
For integers $a$ and $b$ such that $(a,b)=1$ and $0 \leq a \leq b \leq R$ we set 
\begin{equation}\label{major_arc}
\mathfrak{M}(b,a) = \left\{ \alpha \in [0,1] : \left\vert \alpha - \frac{a}{b}\right\vert \leq \frac{R}{n \log\log R} \right\}.
\end{equation}
Let $\mathfrak{m}$ denote the set of real numbers $\alpha$ for which there do not exist integers $a$ and $b$ such that 
$(a,b)=1$, $1 \leq b < R$ , and $\alpha \in \mathfrak{M}(b,a)$.

\begin{lemma}\label{minorestimate} 
For $\alpha \in \mathfrak{m}$ and large $n$,  
\begin{equation}
S_{d,n}(\alpha) \ll \frac{dn}{\phi(d)} \cdot \frac{\log \log R}{R} .
\end{equation}
\end{lemma}

\begin{proof}
This is a restatement of Lemma 9 from \cite{Sarkozy3}.
\end{proof}

\begin{lemma}\label{sumovermajorarc} 
Let $a$ and $b$ be integers such that $0 \leq a \leq b \leq R$ and $(a,b)=1$. 
Then for $n$ sufficiently large
\[
 \sum_{t/n \in  \mathfrak{M}(b,a)} \left|S_{d,n}\left(t/n\right)\right| \ll \frac{dn}{\phi(d)\phi(b)}\log R.
\]
\end{lemma}

\begin{proof}
Suppose that $t/n \in \mathfrak{M}(b,a)$. Then 
\begin{equation*}
	\left|\frac{t}{n}-\frac{a}{b}\right| \leq \frac{R}{n\log\log R} \leq \frac{\log n}{n},
\end{equation*}
and since $b \leq R \leq \log n$ we can, for large enough $n$, apply Lemma~\ref{majorestimates} 
with $\alpha$ replaced by $t/n$.
  
Let $u$ and $v$ be integers such that 
\[
 \frac{u}{n} < \frac{a}{b} < \frac{v}{n}, \quad \; v-u =2.
\]
Applying Lemma~\ref{majorestimates} we obtain
\[
\sum_{\substack{t/n \in  \mathfrak{M}(b,a) \\ u/n \leq t/n \leq v/n}} \left|S_{d,n}\left(t/n\right)\right| \ll \frac{dn}{\phi(d)\phi(b)}.
\]
For $t/n \in \mathfrak{M}(b,a)$ with $t/n < u/n$, Lemma~\ref{majorestimates} implies
\[
 \left|S_{d,n}\left(t/n\right)\right| \ll \frac{d}{\phi(d)\phi(b)}\left| \frac{t}{n} - \frac{a}{b}\right|^{-1}
				\ll \frac{d}{\phi(d)\phi(b)}\left| \frac{t}{n} - \frac{u}{n}\right|^{-1}.
\]
Therefore
\begin{align*}
 \sum_{\substack{t/n \in  \mathfrak{M}(b,a) \\ t/n < u/n}} \left|S_{d,n}\left(t/n\right)\right| 
&\ll \frac{dn}{\phi(d)\phi(b)}\sum_{\substack{t/n \in  \mathfrak{M}(b,a) \\ t/n < u/n}}\frac{1}{|t-u|} \\
&\ll \frac{dn}{\phi(d)\phi(b)}\sum_{1 \leq m \leq R/\log\log R} \frac{1}{m} \ll \frac{dn}{\phi(d)\phi(b)}\log R.
\end{align*}
Similarly
\[
 \sum_{\substack{t/n \in  \mathfrak{M}(b,a) \\ v/n < t/n}} \left|S_{d,n}\left(t/n\right)\right| \ll \frac{dn}{\phi(d)\phi(b)}\log R.
\]
The result follows.
\end{proof}

A multiplicative arithmetic function $f$ is called strongly multiplicative if 
$f(p^k) = f(p)$ for every prime $p$ and positive integer $k$. The next lemma contains a standard deduction
on the average order over arithmetic progressions for certain strongly mutliplicative arithmetic functions.  

\begin{lemma}\label{average_order} 
Let $x$ be a real number such that $x \geq 1$, and let $d$ and $r$ be positive integers.
If $f$ is a strongly multiplicative arithmetic function such that $f(m) \geq 1$ for every positive integer $m$ 
and $f(p) = 1 + O(p^{-1})$. Then 
\[
 \sum_{\substack{m \leq x \\ m \equiv r \mod{d}}} f(m) \ll f((r,d))\frac{x}{d}.
\]
\end{lemma}

\begin{proof} 
Let $g$ be the arithmetic function defined by
\[
 	g(m) = \sum_{k|m}\mu\left(\frac{m}{k}\right)f(k),
\]
where $\mu$ is the M\"{o}bius function.
Using the fact that $f$ is strongly multiplicative we deduce that 
\[
 		g(m) = \mu(m)^2 \prod_{p | m} (f(p) -1). 
\]
Since $f(m) \geq 1$ for every positive integer $m$ it follows that $g$ is a non-negative valued arithmetic function. 
By the M\"{o}bius inversion formula $f(m) = \sum_{k|m}g(k)$, therefore
\begin{equation*}
\sum_{\substack{m \leq x \\ m \equiv r \mod{d} }} f(m) = \sum_{\substack{m \leq x \\ m \equiv r \mod{d}}} \sum_{k |m} g(k)
	= \sum_{k \leq x} g(k) \sum_{\substack{m \leq x \\ m \equiv r \mod{d} \\ m \equiv 0 \mod{k}}} 1.
\end{equation*}
The last sum above is zero if $(k,d) \nmid r$ and at most $x(d,k)/(dk)$ if $(k,d)|r$. This implies, 
since $g$ is a non-negative valued function, that
\begin{align*}
 \sum_{\substack{m \leq x \\ m \equiv r \mod{d} }} f(m)  
&\leq \frac{x}{d}  \sum_{\substack{k \leq x \\ (k,d)|r}} \frac{g(k)(k,d)}{k} 
= \frac{x}{d}\sum_{s|(r,d)} s \sum_{\substack{k \leq x \\ (k,d)=s}} \frac{g(k)}{k} \\
&= \frac{x}{d} \sum_{s|(r,d)} \sum_{\substack{l \leq x/s \\ (l,d/s)=1}} \frac{g(sl)}{l}.
\end{align*}
For positive integers $u$ and $v$ it can be verified that $g(uv) \leq g(u)g(v)$, thus
\begin{align*}
 \sum_{\substack{m \leq x \\ m \equiv r \mod{d} }} f(m)  
	&\leq \frac{x}{d} \sum_{s|(r,d)} g(s)\sum_{l \leq x} \frac{g(l)}{l} \\
	&\leq f((r,d))\frac{x}{d} \prod_{p \leq x} \left(1 + \frac{g(p)}{p}\right) \\
	&=  f((r,d))\frac{x}{d} \prod_{p \leq x} \left(1 + \frac{f(p)-1}{p}\right).
\end{align*}
Since $f(p) \geq 1$ and $f(p) = 1 + O(p^{-1})$ the previous product is bounded from above by the
absolutely convergent infinite product $\prod_{p} (1 + p^{-1}(f(p)-1))$. Therefore
\begin{equation*}
 \sum_{\substack{m \leq x \\ m \equiv r \mod{d} }} f(m)  \ll f((r,d))\frac{x}{d}.
\end{equation*}
\end{proof}

The next lemma is analogous to Proposition 11 of Green \cite{Green}. 

\begin{lemma}\label{sumfourthpowers}
\begin{equation*}
 \sum_{t=0}^{n-1} |S_{d,n}(t/n)|^4 \ll \left(\frac{dn}{\phi(d)}\right)^4.
\end{equation*}
\end{lemma}

\begin{proof}
By Gallagher's inequality \cite[Lemma 1.2]{Montgomery} we have  
\[
 \sum_{t=0}^{n-1} |S_{d,n}(t/n)|^4 \leq n \int_0^1 |S_{d,n}(\alpha)|^4 d\alpha + 2 \int_0^1 |S_{d,n}(\alpha)^3 S_{d,n}'(\alpha)|d\alpha,
\]
where $S_{d,n}'(\alpha)$ is the derivative of $S_{d,n}(\alpha)$ with respect to $\alpha$. By H\"{o}lder's inequality 
\[
 \int_0^1 |S_{d,n}(\alpha)^3 S_{d,n}'(\alpha)|d\alpha \leq \left(\int_0^1 |S_{d,n}(\alpha)|^4 d\alpha\right)^{3/4}
						\left(\int_0^1 |S_{d,n}'(\alpha)|^4 d\alpha\right)^{1/4}.
\]
Let $r_d(m)$ denote the number of pairs $(p_1,p_2)$ where $p_1$ and $p_2$ are primes such that 
$p_1,p_2 \equiv 1 \pmod{d}$ and  
\[
 	\frac{p_1-1}{d} + \frac{p_2-1}{d} = m.
\]
By Parseval's identity,  
\[
\int_0^1 |S_{d,n}(\alpha)|^4 d\alpha \leq (\log n)^4 \sum_{m \leq n} r_d(m)^2 
\]
and 
\[
 \int_0^1 |S_{d,n}'(\alpha)|^4 d\alpha \leq 2\pi (n \log n)^4 \sum_{m \leq n} r_d(m)^2 .
\]
From the above we deduce that
\begin{equation}\label{bridge}
 \sum_{t=0}^{n-1} |S_{d,n}(t/n)|^4 \ll n(\log n)^4 \sum_{m \leq n} r_d(m)^2.
\end{equation}
For each positive integer $m$ we have 
\[
 r_d(m) \leq \big|\{ \; p \; : \; 1 < p \leq dm+2, \; p \equiv 1 \mod{d}, \; \text{$dm+2-p$ is a prime} \; \}\big|.  
\]
To bound $r_d(m)$ we apply the combinatorial sieve to estimate the size of the set above. In particular,
Corollary 2.4.1 of \cite{Halberstam_&_Richert} implies
\[
 	r_d(m) \ll \prod_{p | d(dm+2)} \left(1 - \frac{1}{p}\right)^{-1} \frac{dm+1}{\phi(d)\log^2((dm+1)/d)}.		
\]
Note that 
\[
\prod_{p | d(dm+2)} \left(1 - \frac{1}{p}\right)^{-1} \leq \frac{d}{\phi(d)} \left(\frac{dm+2}{\phi(dm+2)}\right),
\]
therefore
\[
 	r_d(m) \ll \frac{d^2 m}{\phi(d)^2(\log m)^2} \left(\frac{dm+2}{\phi(dm+2)}\right).
\]
This implies 
\[
 \sum_{m \leq n} r_d(m)^2 \ll \frac{d^4n^2}{\phi(d)^4(\log n)^4} 
	\sum_{\substack{u \leq dn+2 \\ u \equiv 2 \mod{d}}} \left(\frac{u}{\phi(u)} \right)^2.
\]
Let $f(u) = (u/\phi(u))^2$. It can verified that $f$ is a strongly multiplicative arithmetic function such that 
$f(u) \geq 1$ for every positive integer $u$ and $f(p) = 1 + O(p^{-1})$. 
Thus, we can apply Lemma~\ref{average_order} to obtain  
\[
 \sum_{\substack{u \leq dn+2 \\ u \equiv 2 \mod{d}}} \left(\frac{u}{\phi(u)} \right)^2 \ll n.
\]
Therefore 
\[
  \sum_{m \leq n} r_d(m)^2 \ll \frac{d^2 n^3}{\phi(d)^2 (\log n)^4},
\]
and thus, on account of (\ref{bridge}), the result follows.
\end{proof}

\section{A density increment}
Throughout this section $n$ denotes a positive integer and $A$ a subset of $\{1, \ldots , n\}$. For any real $\alpha$ we set 
\[
F(\alpha) = \sum_{a \in A} e(\alpha a) , \quad \quad  F_1(\alpha) = \sum_{\substack{a \in A \\ a \leq n/2}} e(\alpha a).
\]
We denote by $C_1$ a fixed positive constant. This constant will be used throughout the rest of 
the paper. We will need $C_1$ to be sufficiently large, but 
it should be noted that the size of $C_1$ will never be determined by $n$ or $A$. Let $\delta$ denote the density of $A$, that 
is, $|A| = \delta n$. The following parameters are defined in terms of $C_1$ and $\delta$.
\begin{equation}\label{definition_of_R}
 R(\delta) = (C_1\delta^{-1})^{(\log\log C_1\delta^{-1})^{7/8}}, 
\end{equation}
\begin{equation}\label{definition_of_theta}
 \theta(\delta) =  (C_1 \delta^{-1})^{-4(\log\log\log C_1 \delta^{-1})^{-1}}. 
\end{equation}
\begin{equation}\label{definition_of_Q1}
Q_1 = (C_1\delta^{-1})^{(\log \log C_1\delta^{-1})^{1/8}},
\end{equation}
\begin{equation}\label{definition_of_Lambda}
		\Lambda = \left[ \frac{3}{4}\log\log\log C_1\delta^{-1} \right],
\end{equation}

With $R=R(\delta)$ we let $\mathfrak{M}(q,a)$ be defined as in (\ref{major_arc}), and 
for any positive integer $q \leq R$ we set 
\[
	\mathfrak{M}(q) = \bigcup_{\substack{a=0 \\ (a,q)=1}}^q \mathfrak{M}(q,a).
\]

\begin{lemma}\label{lemma_nonuniform} Let $d$ be a positive integer such that $d \leq \log n$. Suppose that 
$A-A$ does not intersect $\mathcal{S}_d$ and that 
\begin{equation}\label{size_constraint}
	C_1\delta^{-1} \leq  e^{(\log\log n)^{1/2}}.
\end{equation}
Provided $C_1$ and $n$ are sufficiently large there exists a positive integer $q \leq R(\delta)$ such that  
\begin{equation}\label{arith_nonuniform}
\sum_{\substack{t=1 \\ t/n \in \mathfrak{M}(q)}}^{n-1} \left|F\left(t/n\right)\right|^2  \geq \theta(\delta) |A|^2. 
\end{equation}
\end{lemma}

\begin{proof}
Here we adopt the method used in \cite{BPPS}.
Given any positive integer $\lambda$ we make the following definitions. 
For integers $a$ and $k$, with $k \geq 1$, we define 
\[
\mathfrak{M}_\lambda(k,a) 
= \left\{ \alpha \in [0,1] : \left\vert \alpha - \frac{a}{k}\right\vert \leq \frac{\lambda R}{n \log\log R} \right\},
\]
and for real numbers $K,U \geq 1$ we define
\[
\mathcal{P}_\lambda(K,U) = \left\{ \frac{a}{k} : 1 \leq a \leq k \leq K, (a,k)=1, 
			\max_{t/n \in \mathfrak{M}_\lambda(k,a)}\left| F_1(t/n) \right| \geq |A|/U \right\}.
\]
Furthermore, we set 
\begin{equation}\label{Ql}
		Q_{\lambda} =  Q_1^{2^\lambda -1} 
\end{equation}
and
\[
 \mu_\lambda = \max_{\substack{1 \leq K \leq Q_\lambda \\ 1 \leq U}} \frac{|\mathcal{P}_\lambda(K,U)|}{U^2}. 
\]
Let $K_\lambda$ and $U_\lambda$ denote a pair for which $\mu_\lambda$ takes its maximum.
As $K=U=1$ is considered in the definition of $\mu_\lambda$ we have 
\begin{equation}\label{mu_ineq1}
	1 \leq \mu_{\lambda} \leq \frac{K_\lambda^2}{U_\lambda^2}.
\end{equation}
It follows that 
\begin{equation}\label{chain}
	1 \leq U_{\lambda} \leq K_\lambda \leq Q_\lambda.
\end{equation}

For each $\lambda \leq \Lambda$ we want that the intervals $\mathfrak{M}_\lambda(k,a)$ with $k \leq Q_\lambda$ 
to be pairwise disjoint. It can be verified that this will happen if 
\begin{equation}\label{wish}
	\frac{2\lambda R}{n\log \log R} < \frac{1}{Q_\lambda^2} \quad \quad (\text{for $\lambda \leq \Lambda$}).
\end{equation}
To show this is true we estimate $\lambda$, $R$, and $Q_\lambda$ for $\lambda \leq \Lambda$. 
By (\ref{definition_of_Lambda}) and (\ref{size_constraint}) we deduce that 
\[
		\lambda \leq \frac{3}{4}\log\log\log\log n  \quad \quad (\text{for $\lambda \leq \Lambda$}). 
\]
By (\ref{definition_of_Lambda}) we find that $2^{\lambda} \leq (\log\log C_1 \delta^{-1})^{3/4}$, and thence by $(\ref{definition_of_Q1})$ and $(\ref{Ql})$ we find that 
\[
	\log Q_{\lambda} \leq  2^{\lambda} \log Q_1 \leq (\log \log C_1 \delta^{-1})^{7/8}\log C_1\delta^{-1}.
\] 
By (\ref{definition_of_R}) this implies $\log Q_{\lambda} \leq \log R$, and so
\begin{equation}\label{okay}
 	Q_{\lambda} \leq R.
\end{equation}
By (\ref{definition_of_R}) and (\ref{size_constraint}) we find, for $n$ large enough, that
\begin{equation}\label{R_home}
 3 \leq R \leq \log n.
\end{equation}
From the above estimates for $\lambda$, $R$, and $Q_\lambda$ we deduce that (\ref{wish}) holds 
for sufficiently large $n$. Therefore, when $\lambda \leq \Lambda$ we have   
\[
		\mu_{\lambda}|A|^2 = |\mathcal{P}_\lambda(K_\lambda,U_\lambda)| \frac{|A|^2}{U_\lambda^2} 
		\leq \sum_{t=0}^{N-1} \left|F_1(t/n)\right|^2 \leq n|A|.
\]	
So
\begin{equation}\label{mu_ineq2}
			\delta \leq \mu_{\lambda}^{-1}.
\end{equation}

Let us assume, to obtain a contradiction, that 
\begin{equation}\label{assumption}
\sum_{\substack{t=1 \\ t/n \in \mathfrak{M}(q)}}^{n-1} \left|F(t/n)\right|^2  < \theta(\delta) |A|^2  \quad 
(\text{for all $1 \leq q \leq R$}). 
\end{equation}
By using Lemma~\ref{lemma_CR} and (\ref{assumption}) we will show, provided $C_1$ and $n$ are sufficiently large, that 
\begin{equation}\label{target}
	\mu_{\lambda+1} \geq \theta(\delta)^{-1/2} \mu_{\lambda} \quad \text{(for $1 \leq \lambda \leq \Lambda$)}.
\end{equation}

Assuming for now that (\ref{target}) holds we show how a contradiction is obtained, thus proving that the assumption 
(\ref{assumption}) is false.   
Since $\mu_1 \geq 1$, it follows from (\ref{target}) that $\mu_{\Lambda+1} \geq \theta(\delta)^{-(1/2)\Lambda}$, and thus by (\ref{mu_ineq2}) we have 
\[
				\delta \leq \theta(\delta)^{(1/2)\Lambda}.
\]
We can take $C_1$ to be large enough so that (\ref{definition_of_Lambda}) implies $\Lambda \geq (1/4)\log_3 C_1 \delta^{-1}$, 
then by (\ref{definition_of_theta}) we find that
\[
\delta \leq C_1^{-1}\delta < \delta,
\]
a contradiction. Therefore (\ref{assumption}) cannot hold for all $1 \leq q \leq R$. 

We now proceed to show that (\ref{target}) holds.  To that end, let us fix $\lambda$ with $1 \leq \lambda \leq \Lambda$. 
For now we also fix a rational $a/k$ in $\mathcal{P}_\lambda(U_\lambda,K_\lambda)$. We associate with $a/k$ a fraction 
$u/n \in \mathfrak{M}_\lambda(k,a)$ such that $|F(u/n)| \geq |A|/U_\lambda$. Such a $u/n$ exists by the way $a/k$ was chosen. 

Since $A-A$ contains no integers from $\mathcal{S}_d$ we find that 
\[
	\sum_{t = 0}^{n-1} F_1(u/n+ t/n)F(-t/n)S_{d,n}(t/n) = 0.
\]
By the triangle inequality, Lemma~\ref{Siegel_Walfisz}, and the way $u/n$ was chosen we find that
\begin{equation}\label{start}
\frac{|A|^2}{U_\lambda} \cdot \left(\frac{d n}{\phi(d)} \right) \ll  
                           \sum_{t=1}^{n-1} |F_1(u/n+t/n)||F(t/n)||S_{d,n}(t/n)|.
\end{equation}

Set 
\begin{equation}\label{definition_of_Y}
 	Y = (C_1\delta^{-1})^{3/2} Q_\lambda^2 
\end{equation}
and let $\mathcal{N}$ denote the set of $t/n$ such that $|F(t/n)| \leq |A|/Y$.
By two applications of the Cauchy-Schwartz inequality, Parseval's identity, and 
Lemma~\ref{sumfourthpowers} we find that
\begin{align*}
  & \sum_{t/n \in \mathcal{N}} |F_1(u/n+t/n)||F(t/n)||S_{d,n}(t/n)| \\
&\leq \left( \sum_{t=0}^{n-1} |F_1(u/n + t/n)|^2 \right)^{1/2} 
				\left( \sum_{t/n \in \mathcal{N}}|F(t/n)|^4 \right)^{1/4} 
				\left( \sum_{t=0}^{n-1} |S_{d,n}(t/n)|^4 \right)^{1/4} \\
&\ll \frac{dn^{3/2}|A|^{1/2}}{\phi(d)}\left( \sum_{t/n \in \mathcal{N}}|F(t/n)|^4 \right)^{1/4}.
\end{align*}
Now  
\begin{align*}
\left( \sum_{t/n \in \mathcal{N}}|F(t/n)|^4 \right)^{1/4} 
&\leq \max_{t/n \in \mathcal{N}} |F(t/n)|^{1/2} 
					\left(\sum_{t=0}^{n-1} |F(t/n)|^2 \right)^{1/4} \\
&\leq \frac{|A|^{1/2}}{Y^{1/2}} (n|A|)^{1/4} = \frac{n^{1/4}|A|^{3/4}}{Y^{1/2}}.
\end{align*}
Therefore 
\[
\sum_{t/n \in \mathcal{N}} |F_1(u/n+t/n)||F(t/n)||S_{d,n}(t/n)| \ll \frac{dn^{7/4}|A|^{5/4}}{\phi(d) Y^{1/2}}.
\]
By (\ref{chain}) and (\ref{definition_of_Y}) we find that 
\[
 Y^{-1/2} = C_1^{-3/4}\delta^{3/4}Q_\lambda^{-1} \leq C_1^{-3/4}|A|^{3/4}n^{-3/4}U_\lambda^{-1},
\]
thus
\begin{equation}\label{small_contribution_0}
\sum_{t/n \in \mathcal{N}} |F_1(u/n+t/n)||F(t/n)||S_{d,n}(t/n)| 
\ll C_1^{-3/4}\frac{|A|^2}{U_\lambda}\left(\frac{dn}{\phi(d)}\right).
\end{equation}

Let $\mathcal{N}_1$ denote the set of $t/n$ such that $|F_1(u/n+ t/n)| \leq |A|/Y$.  
By the same reasoning used in the deduction of (\ref{small_contribution_0}) we find that 
\begin{equation}\label{small_contribution_1}
\sum_{t/n \in \mathcal{N}_1} |F_1(u/n+t/n)||F(t/n)||S_{d,n}(t/n)| 
\ll C_1^{-3/4}\frac{|A|^2}{U_\lambda}\left(\frac{dn}{\phi(d)}\right).
\end{equation}

For $\lambda \leq \Lambda$ we have $Q_{\lambda+1}/Q_\lambda < R$. Indeed, (\ref{definition_of_Lambda}) 
and (\ref{Ql}) imply
\[
	\frac{Q_{\lambda+1}}{Q_\lambda} \leq Q_1^{2^\Lambda} 
						\leq (C_1\delta^{-1})^{(\log\log C_1\delta^{-1})^{3/4}} < R.
\]
Let $\mathfrak{m}^\ast$ denote the union of the $\mathfrak{M}(q)$ with $Q_{\lambda+1}/Q_{\lambda} \leq q \leq R$.
By the Cauchy-Schwartz inequality we find that 
\begin{equation}\label{minortrick1}
\sum_{t/n \in \mathfrak{m}^\ast} |F_1(u/n+t/n)||F(t/n)||S_{d,n}(t/n)| 
		\leq (n|A|)\sup_{t/n \in \mathfrak{m}^\ast_\lambda}|S_{d,n}(t/n)|.
\end{equation}
We are now going to show that 
\begin{equation}\label{super}
 \sup_{t/n \in \mathfrak{m}^\ast_\lambda}|S_{d,n}(t/n)|  \ll 
 					C_1^{-1} U_\lambda^{-1} \delta \left(\frac{dn}{\phi(d)}\right).
\end{equation}
Suppose that $t/n \in \mathfrak{m}^\ast$, then $t/n \in \mathfrak{M}(q,a)$ for 
some integers $a$ and $q$ such that $0 \leq a \leq q$, $(a,q)=1$, and $Q_{\lambda+1}/Q_{\lambda} \leq q \leq R$. 
Since $q \leq R \leq \log n$, we deduce from Lemma~\ref{majorestimates} that 
\[
             S_{d,n}(t/n) \ll \frac{dn}{\phi(d)\phi(q)}. 
\]
Using the well-known estimate
\begin{equation}\label{phi_estimate}
 	\phi(q) \gg \frac{q}{\log\log q},
\end{equation}
(see for example \cite[Theorem 328]{Hardy_&_Wright}), we obtain
\begin{equation}\label{v6}
             S_{d,n}(t/n) \ll  \left( \frac{dn}{\phi(d)} \right)\frac{\log\log q}{q}.
\end{equation}
The lower bound on $q$ implies
\begin{equation}\label{v7}
\frac{\log\log q}{q} \ll  \frac{\log\log Q_{\lambda+1}/Q_{\lambda}}{Q_{\lambda+1}/Q_{\lambda}}.
\end{equation}
By (\ref{Ql}) we have $Q_{\lambda+1}/Q_{\lambda} = Q_{\lambda}Q_1 = Q_1^{2^\lambda}$, thus
\[
 \frac{\log \log Q_{\lambda+1}/Q_{\lambda}}{Q_{\lambda+1}/Q_{\lambda}} 
 = \frac{\log \log Q_1^{2^\lambda}}{Q_{\lambda}Q_1}
 		= \frac{\lambda (\log 2) + \log \log Q_1}{Q_\lambda Q_1} .
\]
Using (\ref{definition_of_Q1}) and (\ref{definition_of_Lambda}) we find that $\lambda \ll \log \log Q_1$, by this and (\ref{chain}) we obtain
\[
\frac{\log \log Q_{\lambda+1}/Q_{\lambda}}{Q_{\lambda+1}/Q_{\lambda}} \ll \frac{\log \log Q_1}{U_\lambda Q_1}.
\]
Using (\ref{definition_of_Q1}) we find, by taking $C_1$ large enough, that 
\[
 	\log \left( \frac{\log \log Q_1}{Q_1} \right) \leq -\log C_1\delta^{-1},
\] 
and thus
\[
\frac{\log \log Q_1}{Q_1} \leq C_1^{-1} \delta.
\]
From (\ref{v7}) and the subsequent estimates we obtain 
\begin{equation}\label{v11}
 \frac{\log\log q}{q} \ll C_1^{-1}U_\lambda^{-1}\delta,
\end{equation}
Since $t/n \in \mathfrak{m}^\ast$ is arbitrary (\ref{v6}) and (\ref{v11})
imply that (\ref{super}) is true. By (\ref{minortrick1}) and (\ref{super}) we have 
\begin{equation}\label{minorstar_contribution}
\sum_{t/n \in \mathfrak{m}^\ast} |F_1(u/n+t/n)||F(t/n)||S_{d,n}(t/n)| 
\ll C_1^{-1}\frac{|A|^2}{U_\lambda}\left(\frac{dn}{\phi(d)}\right).
\end{equation}

The contribution to the sum in (\ref{start}) coming from the terms with $t/n \in \mathfrak{m}$ 
can similarly be bounded. By the Cauchy-Schwartz inequality and Lemma~\ref{minorestimate} we find that 
\begin{align*}
\sum_{t/n \in \mathfrak{m}} |F_1(u/n+t/n)||F(t/n)||S_{d,n}(t/n)| 
&\leq (n|A|)\sup_{t/n \in \mathfrak{m}}|S(t/n)| \\
&\ll (n|A|)\left( \frac{dn}{\phi(d)} \right) \frac{\log\log R}{R}.
\end{align*}
Since $R \geq Q_{\lambda+1}/Q_\lambda$ the argument used the previous paragraph implies 
\begin{equation}\label{minor_contribution}
\sum_{t/n \in \mathfrak{m}} |F_1(u/n+t/n)||F(t/n)||S_{d,n}(t/n)| 
\ll C_1^{-1}\frac{|A|^2}{U_\lambda}\left(\frac{dn}{\phi(d)}\right).
\end{equation}

Let $\mathfrak{N}(b,a)$ be the set of $t/n \in \mathfrak{M}(b,a)$ with $t/n \ne 0$ such that 
\[
|F(t/n)| \geq \frac{|A|}{Y}, \quad \quad |F_1(u/n+t/n)| \geq \frac{|A|}{Y}.
\]
By (\ref{small_contribution_0}), (\ref{small_contribution_1}),  (\ref{minorstar_contribution}), and 
(\ref{minor_contribution}) it follows for $C_1$ large enough that
\begin{align*}
& \frac{d|A|^2 n}{\phi(d)U_\lambda} \ll \\ 
& \sum_{b \leq Q_{\lambda+1}/Q_\lambda} \sum_{(a,b)=1} \max_{t/n \in \mathfrak{N}(b,a)} |F(t/n)| 
             \max_{t/n \in \mathfrak{N}(b,a)} |F_1(u/n+t/n)| \sum_{t/n \in \mathfrak{M}(b,a)} |S_{d,n}(t/n)|.
\end{align*}
Since $d \leq \log n$ we can apply Lemma~\ref{sumovermajorarc} to the inner sum above to obtain 
\[
\frac{|A|^2 }{U_\lambda \log R} \ll \sum_{b \leq Q_{\lambda+1}/Q_\lambda} \frac{1}{\phi(b)}
\sum_{(a,b)=1}  \max_{t/n \in \mathfrak{N}(b,a)} |F(t/n)|\max_{t/n \in \mathfrak{N}(b,a)} |F_1(u/n+t/n)|.
\]

Let $\mathcal{L}(L,V,W)$ denote the set of reduced fractions $b/l \in [0,1]$ such that 
\[
							\frac{L}{2} \leq l \leq L,
\]
\[
\frac{|A|}{V} \leq \max_{t/n \in \mathfrak{M}(l,b)}|F(t/n)| \leq 2\frac{|A|}{V},
\]
\[
\frac{|A|}{W} \leq \max_{t/n \in \mathfrak{M}(l,b)}|F_1(u/n+t/n)| \leq 2\frac{|A|}{W}.
\]
For $b/l \in \mathcal{L}(L,V,W)$, we have  
\[
 \frac{1}{\phi(l)} \max_{t/n \in \mathfrak{M}(l,b)} |F(t/n)|\max_{t/n \in \mathfrak{M}(l,b)} |F_1(u/n+t/n)|
 	\ll \frac{(\log\log 3L)|A|^2}{LVW}
\]
by (\ref{phi_estimate}). Therefore
\[
\frac{|A|^2 }{U_\lambda \log R} 
\ll \sum_{L}\sum_{V} \sum_{W}|\mathcal{L}(L,V,W)|\frac{(\log\log 3L)|A|^2}{LVW}.
\]
where $L$ runs through all the powers of $2$ in the interval $[1,2Q_{\lambda+1}/Q_\lambda]$, 
and $V$ and $W$ run through all the powers of $2$ in the interval $[1,2Y]$.
There must exist a triple $(L,V,W)$ of such indices such that 
\[
	|\mathcal{L}(L,V,W)| 
	\gg  \frac{LVW}{U_\lambda(\log\log 3L)(\log R)}.
\]
We associate this triple with $a/k$.

The number of possible triples $(L,V,W)$ is $\ll \log (Q_{\lambda+1}/Q_\lambda) (\log Y)^2$, which by (\ref{okay}) 
and (\ref{definition_of_Y}) is $\ll (\log R)^3$.
Therefore there exists a subset $\mathcal{K} \subset \mathcal{P}_\lambda$, satisfying  
\begin{equation}\label{size1}
			 |\mathcal{K}| \gg \frac{|P_\lambda(K_\lambda,U_\lambda)|}{(\log R)^3},
\end{equation}
such that for each $a/k \in \mathcal{K}$ we associate the same triple, say $(L,V,W)$. 

Let $a/k \in \mathcal{K}$, then together with the associated fraction $u/n \in \mathfrak{M}_\lambda(k,a)$, 
we associate a set $\mathcal{L}_{a/k}$ of rationals $b/l$, $0 \leq b \leq l$, $(b,l)=1$, $L/2 \leq l \leq L$, such that 
\begin{equation}\label{dm1}
	|\mathcal{L}_{a/k}| \gg \frac{LVW}{U_\lambda (\log \log 3L )(\log R)},
\end{equation}
\begin{equation}\label{dm2}
\frac{|A|}{V} \leq \max_{v/n \in \mathfrak{M}(l,b)} |F(v/n)|       \leq \frac{2|A|}{V},
\end{equation}
\begin{equation}\label{dm3}
\frac{|A|}{W} \leq \max_{w/n \in \mathfrak{M}(l,b)} |F_1(u/n+w/n)| \leq \frac{2|A|}{W}.
\end{equation}
Set 
\[
\mathcal{Q} = \left\{ \; \frac{a}{k}+\frac{b}{l} \; : \; \frac{a}{k} \in \mathcal{K}, \;  \frac{b}{l} \in \mathcal{L}_{a/k} \; \right\}.
\]

Let us estimate the cardinality of $\mathcal{Q}$.  
Since $L \leq Q_{\lambda+1}/Q_\lambda \leq R$, assumption (\ref{assumption}) and (\ref{dm2}) imply
\[
\left|\left\{ \; b \; : \; \frac{b}{l} \in \bigcup \mathcal{L}_{a/k} \; \right\} \right| \left(\frac{|A|}{V}\right)^2 
\leq \sum_{t/n \in \mathfrak{M}(l)}|F(t/n)|^2 \leq \theta(\delta)|A|^2.
\]
So that
\[
 \left|\left\{ \; b \; : \; \frac{b}{l} \in \bigcup \mathcal{L}_{a/k} \; \right\} \right| \ll \theta(\delta) V^2 .
\]
Lemma~\ref{lemma_CR} then implies 
\[
   |\mathcal{Q}| \gg |\mathcal{K}| \cdot \frac{L^2V^2W^2}{U_\lambda^2 (\log \log 3L)^2 (\log R)^2} \cdot  
   			\frac{\theta(\delta)^{-1}}{L V^2 \tau^8 (1+\log K_\lambda)}.
\]
From (\ref{chain}) and (\ref{okay}) we obtain $\log K_\lambda \leq \log R$, by this and (\ref{size1}) it follows that
\begin{equation}\label{road}
	|\mathcal{Q}| \gg  W^2 \left( \frac{\theta(\delta)^{-1}}{\tau^8(\log R)^6 } \right)
						\frac{|\mathcal{P}_\lambda(K_\lambda,U_\lambda)|}{U_\lambda^2}.
\end{equation}

Note that $\mathcal{Q}$ is a subset of $(0,2]$. Let $\mathcal{Q}_1 = \mathcal{Q} \cap (0,1]$ and 
$\mathcal{Q}_2 = \mathcal{Q} \cap (1,2]$. Let us assume without loss of generality that 
$|\mathcal{Q}_1| \geq (1/2)|\mathcal{Q}|$. 
If this is not the case, then $|\mathcal{Q}_2| \geq (1/2)|\mathcal{Q}|$, and we 
can replace $\mathcal{Q}_1$ in the argument below by the rational numbers in $\mathcal{Q}_2$ shifted to the left by $1$.
Since $|\mathcal{Q}_1| \geq (1/2)|\mathcal{Q}|$ we see that (\ref{road}) is still valid with $\mathcal{Q}$ replaced by $\mathcal{Q}_1$ 

Let $r/s = a/k + b/l$ be in $\mathcal{Q}_1$. For $u/n \in \mathfrak{M}_\lambda(k,a)$ and $w/n \in \mathfrak{M}(l,b)$ we have 
\[	
\left| \frac{r}{s} - \left( \frac{u}{n} + \frac{w}{n} \right) \right| \leq
 \left| \frac{u}{n} - \frac{a}{k} \right| + \left|\frac{w}{n}- \frac{b}{l} \right| \leq \frac{(\lambda+1)R}{n\log\log R},
\]
and therefore $u/n + w/n \in  \mathfrak{M}_{\lambda+1}(s,r)$.  Thus, by (\ref{dm3}) we deduce that
\begin{equation}\label{trans}
		\max_{t/n \in \mathfrak{M}_{\lambda+1}(s,r)} |F_1(t/n)| \geq \frac{|A|}{W}  \quad (\text{for $r/s \in \mathcal{Q}_1$}).
\end{equation}
We now estimate the size of the denominator of $r/s$. Certainly $s \leq kl \leq K_\lambda L$. By (\ref{chain}) we have 
$K_\lambda \leq Q_\lambda$ and $L$ was chosen to satisfy $L \leq Q_{\lambda+1}/Q_{\lambda}$. Therefore 
$s \leq Q_{\lambda +1}$ whenever $r/s \in \mathcal{Q}_1$. By this and (\ref{trans}) we obtain 
\begin{equation}\label{inside}
			\mathcal{Q}_1 \subset \mathcal{P}_{\lambda+1}(Q_{\lambda+1},W).
\end{equation}

 By (\ref{road}), with $\mathcal{Q}$ replaced by $\mathcal{Q}_1$, and (\ref{inside}) we find that 
\[
 \frac{|\mathcal{P}_{\lambda+1}(Q_{\lambda+1},W)|}{W^2}| \gg  \left( \frac{\theta(\delta)^{-1}}{\tau^8(\log R)^6 } \right)
						\frac{|\mathcal{P}_\lambda(K_\lambda,U_\lambda)|}{U_\lambda^2}.
\]
This implies
\begin{equation}\label{inc1}
		\mu_{\lambda+1} \gg \frac{\theta(\delta)^{-1}}{\tau^8 (\log R)^6} \mu_\lambda.
\end{equation}

We now estimate $\tau$ the maximum of the divisor function up to $K_\lambda L \leq Q_{\lambda+1}$.
If $d(m)$ is the number of divisors of $m$ then  
\[
 	\log d(m) \ll \frac{\log m}{\log\log m},
\]
(see \cite[Theorem 317]{Hardy_&_Wright}). Thus, by (\ref{Ql}), we have
\[
	\log \tau \ll \frac{\log Q_{\lambda+1}}{\log \log Q_{\lambda+1}} 
		  \ll \frac{2^\lambda \log Q_1}{\log \log Q_1},
\]
and since $\lambda \leq \Lambda$ we deduce from (\ref{definition_of_Q1}) and (\ref{definition_of_Lambda}) that 
\[
 \log \tau \ll \frac{\log C_1\delta^{-1}}{(\log\log C_1\delta^{-1})^{1/4}}.
\]
It follows from (\ref{definition_of_theta}) that 
\begin{equation}\label{appx1}
	 	\log \tau = o(\log \theta(\delta)^{-1}) \quad \quad \text{(for $C_1\delta^{-1} \to \infty$)}. 
\end{equation}
We also find from (\ref{definition_of_R}) and (\ref{definition_of_theta}) that 
\begin{equation}\label{appx2}
		\log \log R = o(\log \theta(\delta)^{-1}) \quad \quad \text{(for $C_1\delta^{-1} \to \infty$)}.  
\end{equation}
Since $\theta(\delta)^{-1}$ tends to infinity as $C_1\delta^{-1}$ tends to infinity, we deduce from 
(\ref{inc1}), (\ref{appx1}), and (\ref{appx2}) that for $C_1$ sufficiently large  
\[
 \mu_{\lambda+1} \geq \theta(\delta)^{-1/2} \mu_\lambda.
\]
Since $\lambda \leq \Lambda$ was arbitrary (\ref{target}) is true, and as shown earlier the lemma 
can be deduced from this.
\end{proof}
We now derive a density increment argument that will be iterated in the next section to prove our theorem.  

\begin{lemma}\label{crank} Let $d$ be a positive integer such that $d \leq \log n$. Suppose that 
$A-A$ does not intersect $\mathcal{S}_d$ and that $\delta$, the density of $A$, satisfies (\ref{size_constraint}).
Provided $C_1$ and $n$ are sufficiently large there exist positive integers $d'$ and $n'$, and a subset $A'$ of $\{1, \ldots ,n'\}$ 
of size $\delta' n'$, such that $A'-A'$ does not intersect $\mathcal{S}_{d'}$, and moreover; 
\[
	d \leq d' \leq R(\delta)d, \; \quad \; \quad  
	R(\delta)^{-2} n \leq n' \leq n,
\]
\[
 			\delta' \geq \delta\big( 1 + 8^{-1} \theta(\delta) \big).
\]
\end{lemma}

\begin{proof} By the hypotheses Lemma~\ref{lemma_nonuniform} implies there exists a positive integer $q \leq R(\delta)$ 
such that (\ref{arith_nonuniform}) is true. With this $q$ and $U= R(\delta)/\log\log R(\delta)$
let $E$ be defined as in Lemma~\ref{bump}. Note that $\mathfrak{M}(q) \subset E$. The inequality (\ref{R_home}) 
is still valid, thus $2\pi q U \leq 2\pi R(\delta)^2 \leq n$ for sufficiently large $n$. 
Therefore, we can apply Lemma~\ref{bump} with $\theta = \theta(\delta)$ to deduce that there exists an arithmetic progression 
$P$ with difference $q$ such that 
\begin{equation}\label{lu1}
				|P| \geq \frac{n \log \log R(\delta)}{32 \pi q R(\delta)}
\end{equation}
and 
\begin{equation}\label{lu2}
		|A \cap P| \geq |P|\delta \big(1+ 8^{-1}\theta(\delta)\big).
\end{equation}
Let $n' = |P|$. Then there exists an integer $c$ and subset $A'$ of $\{1, \ldots , n'\}$ such that  
$A \cap P = \{ \, c+qa' \, : \, a' \in A' \, \}$. Put $d' = dq$. 
Since $A-A$ does not intersect $\mathcal{S}_d$, we deduce that $A'$ does not intersect $\mathcal{S}_{dq}$.
Let the size of $A'$ be  $\delta' n'$. Then (\ref{lu2}) implies
\[
			\delta' \geq \delta\big(1+ 8^{-1}\theta(\delta)\big).
\]
To finish we need to estimate $n'$ and $d'$. Since $q \leq R(\delta)$ we find by (\ref{lu1}) and for $C_1$ large enough 
that $n' \geq R(\delta)^{-2} n$, and clearly, $n' \leq n$. Now, again by the fact that $q \leq R(\delta)$, we 
obtain $q \leq d' = dq \leq R(\delta)q$. This completes the proof.
\end{proof}

\section{Proof of the Theorem}\label{proof} 
Let us assume, for a contradiction, that the theorem is false. Then for $C_1$ and $n$ sufficiently large, there exists 
a subset $A$ of $\{1, \ldots , n\}$ of size $\delta n$, such that $A-A$ does not intersect $\mathcal{S}$ and  
\begin{equation}\label{hyp}
 	\delta  \geq C_1 \left(\frac{\log_2 n}{(\log_3 n)^2(\log_4 n)}\right)^{-\log_5 n}.
\end{equation}

Set
\begin{equation}\label{definition_of_Z}
 Z = \big[ 64 \, \theta(\delta)^{-1}\log C_1\delta^{-1} \big],
\end{equation}
and put $d_0 = 1$, $n_0 = n$, $A_0 = A$, and $\delta_0 = \delta$. 
By using Lemma~\ref{crank} repeatedly we can show that for each integer $k$, with $1 \leq k \leq Z$,
there are integers $d_k$ and $n_k$ and a subset $A_k$ of $\{ 1 , \ldots , n_k \}$ of size $\delta_k n_k$
such that $A_k -A_k$ does not intersect $\mathcal{S}_{d_k}$. Moreover, $d_k$, $n_k$, and $\delta_k$ satisfy
\[
	d_{k-1} \leq d_{k} \leq R(\delta_{k-1})d_{k-1}, \; \quad \; \quad  
	R(\delta_{k-1})^{-2} n_{k-1} \leq n_{k} \leq n_{k-1},
\]
\[
 			\delta_k \geq \delta_{k-1}\big( 1 + 8^{-1} \theta(\delta_{k-1}) \big).
\]
Since $d_0 =1$ and $n_0=n$, these estimates imply
\begin{equation}\label{powertrio}
 d_k \leq R(\delta)^{k}, \quad n_k  \geq R(\delta)^{-2k} n, \quad 
\delta_k \geq \delta\big( 1 + 8^{-1} \theta(\delta) \big)^k.
\end{equation}

Let us show that we can actually perform this iteration $Z$ many times. Let $0 \leq l \leq Z-1$, and suppose 
that we have performed this iteration $l$ many times. To show that 
Lemma~\ref{crank} can be applied a $(l+1)$-th time we need to show that 
$n_l$ is sufficiently large, $d_l \leq \log n_l$, and that (\ref{size_constraint}) is satisfied with $\delta$ replaced by $\delta_l$. 

We begin by estimating $n_l$. By (\ref{powertrio}) we obtain 
\begin{equation}\label{estimate_of_nk}
	\log n_l \geq \log n - 2l \log R(\delta). 
\end{equation}
Since $l < Z$, (\ref{definition_of_R}) and (\ref{definition_of_Z}) imply
\[
l \log R(\delta) \leq 64 \, \theta(\delta)^{-1}(\log C_1\delta^{-1})^2(\log_2 C_1\delta^{-1})^{7/8}.
\]
By (\ref{hyp}) we obtain 
\[
	(\log C_1\delta^{-1})^2(\log_2 C_1\delta^{-1})^{3/4} \leq 2(\log_3 n)^2 (\log_4 n)^{7/8} (\log_5 n)^2
\]
for large enough $n$. By (\ref{definition_of_theta}) and (\ref{hyp}) we find, for $n$ and $C_1$ sufficiently large, that
\[
 \log \theta(\delta)^{-1} = \frac{4\log C_1\delta^{-1}}{\log_3 C_1 \delta^{-1}} 
 \leq \log \left( \frac{\log_2 n}{(\log_3 n)^2(\log_4 n)} \right).
\]
(Here we used that $(\log x)(\log_3 x)^{-1}$ is eventually increasing.) Therefore
\[
\theta(\delta)^{-1} \leq   \frac{\log_2 n}{(\log_3 n)^2(\log_4 n)}. 
\]
From the above we deduce, for $n$ and $C_1$ large enough, that 
\begin{equation}\label{peq2}
	l \log R(\delta) \leq \log_2 n. 
\end{equation}
Therefore, by (\ref{estimate_of_nk}), 
\[
	\log n_l \geq \log n - 2\log_2 n =  \log \left(\frac{n}{(\log n)^2}\right),
\]
and so 
\begin{equation}\label{peq3}
		n_l \geq \frac{n}{(\log n)^2} 
\end{equation}
for $l < Z$. This shows that by taking $n$ to be arbitrarily large, the same is true for $n_l$.

We now show that $d_l \leq \log n_l$. By (\ref{powertrio}) we have $\log d_l \leq l \log R(\delta)$, and 
thus by (\ref{peq2}) we obtain $\log d_l  \leq (1/2)\log_2 n$. For large $n$ this implies
\[
	d_l \leq (\log n)^{1/2} \leq \log \frac{n}{(\log n)^{2}} \leq \log n_l 
\]
by (\ref{peq3}). 

We leave it to the reader to verify that (\ref{hyp}) and (\ref{peq3}) imply, for $n$ and $C_1$ sufficiently large,
that (\ref{size_constraint}) is satisfied with $\delta$ and $n$ replaced by $\delta_l$ and $n_l$ respectively.
Finally, since $A_l - A_l$ does not intersect $\mathcal{S}_{d_l}$ we can apply Lemma~\ref{crank} to obtain the desired outcome.

Since (\ref{powertrio}) is true with $k=Z$ we find that 
\[
 \log \delta_Z \geq Z\log \Big(1+ 8^{-1}\theta(\delta)\Big) - \log C_1\delta^{-1}.
\]
Since $8^{-1}\theta(\delta) < 1$, this implies
\begin{equation}\label{nearly}
 \log \delta_Z \geq  {16}^{-1}Z \theta(\delta) - \log C_1\delta^{-1}.
\end{equation}
(Here we used $\log(1+x) \geq x/2$ for $0 \leq x \leq 1$.) 
For $C_1$ large enough $Z \geq 32 \theta(\delta)^{-1} \log C_1 \delta^{-1}$, thus 
\[
 		\log \delta_Z \geq 2\log C_1\delta^{-1} - \log C_1 \delta^{-1} > 0.
\]
This implies $\delta_Z > 1$, a contradiction, since by definition $\delta_Z \leq 1$. This contradiction establishes the theorem.

\section*{acknowledgements}
The author was supported by a postdoctoral fellowship from the 
Centre de recherches math\'{e}matiques at Montr\'{e}al.

\bigskip
\begin{scriptsize}
\noindent Centre de recherches math\'{e}matiques\\
Universit\'{e} de Montr\'{e}al \\
Case postale 6128, Succursale Centre-ville \\
Montr\'{e}al, H3C 3J7 \\
Canada\\
\end{scriptsize}
\end{document}